\documentclass[12pt]{amsart}
\usepackage{amsmath,amssymb}
\usepackage[dvips]{graphicx}
\newcommand{\N}{\mathbb{N}}
\newcommand{\Z}{\mathbb{Z}}
\newcommand{\R}{\mathbb{R}}

\newcommand{\CU}{\mathcal{U}}
\newcommand{\BF}{\mathbf{F}}

\newcommand{\hooklongrightarrow}{\lhook\joinrel\to}

\DeclareMathOperator{\Int}{Int}

\DeclareMathOperator{\st}{st}
\DeclareMathOperator{\Sd}{Sd}
\DeclareMathOperator{\diam}{diam}
\newtheorem{Th}{Theorem}[section]
\newtheorem{Def}[Th]{Definition}
\newtheorem{Prop}[Th]{Proposition}
\newtheorem{Cor}[Th]{Corollary}
\newtheorem{L'a}[Th]{Lemma}

\newtheorem{Ex}[Th]{Example}
\newtheorem{Claim}{Claim}
\newtheorem{Q}[Th]{Question}
\begin{document}
\title[Non-separable Hilbert manifolds of continuous mappings]
{Non-separable Hilbert manifolds of continuous mappings}
\author[A. Yamashita]{Atsushi Yamashita}
\address{Graduate School of Mathematical Sciences, The University of Tokyo,
3--8--1, Komaba, Meguro-ku, Tokyo 153--8914, Japan} 
\email[A. Yamashita]{yonster@ms.u-tokyo.ac.jp}

\begin{abstract}
Let $X, Y$ be separable metrizable spaces, where $X$ is noncompact
and $Y$ is equipped with an admissible complete metric $d$.
We show that the space $C(X,Y)$ of continuous maps from $X$ into $Y$
equipped with the uniform topology is locally homeomorphic
to the Hilbert space of weight $2^{\aleph_0}$
if (1) $(Y, d)$ is an ANRU, a uniform
version of ANR and (2) the diameters of components of $Y$ is bounded
away from zero. The same conclusion holds for the subspace $C_B(X,Y)$
of bounded maps if $Y$ is a connected complete Riemannian manifold.
\end{abstract}
\maketitle
\section{Introduction}
Most function spaces are infinite-dimensional by nature,
and it is natural to ask whether their local structure is
similar to that of Hilbert spaces or Banach spaces.
From the topological viewpoint, it is known that
two Fr\'{e}chet spaces (i.e., completely metrizable, 
locally convex real topological vector spaces) are homeomorphic
if they have the same weight, due to the characterization 
of Hilbert spaces by Toru\'{n}czyk~\cite{Torunczyk} (cf.~\cite{Torunczyk2}).
Therefore, any Fr\'{e}chet space of weight $\tau$
is homeomorphic to the \textit{Hilbert space} $\ell^2(\tau)$ of weight $\tau$,
in other words, to the Hilbert space with $\tau$ orthonormal basis vectors.

Thus we may expect that many function spaces are
locally homeomorphic to the Hilbert space $\ell^2(\tau)$ for a suitable $\tau$.
Such spaces are called topological \textit{Hilbert manifolds}, or
more precisely, \textit{$\ell^2(\tau)$-manifolds}.

Using the Toru\'{n}czyk's characterization mentioned above,
many function spaces are known to be
Hilbert manifolds. One typical example is the following:
the function space $C(X, Y)$ of continuous maps
from an infinite compact metric space $X$
into a separable complete metric ANR $Y$
with no isolated points is an $\ell^2$-manifold,
with respect to the uniform topology(Sakai~\cite{Sakai}).\footnote{In
this case, the uniform topology coincides with many other
interesting topologies including the compact-open topology,
since $X$ is compact.}

In this paper, we prove corresponding results
for $C(X,Y)$ with the uniform topology
when $X$ is a \textit{non}compact space,
with suitable assumption added to $Y$.
The most significant change from the previous setting is that
the function space $C(X,Y)$ becomes non-separable.\footnote{
In this case, the compact-open topology is coarser than the uniform
topology, and with this topology, $C(X,Y)$ need not have nice
local behavior. This problem is discussed in~\cite{Smrekar-Yamashita}.}

When we apply Toru\'{n}czyk's characterization to prove
a space to be a Hilbert manifold, we have first to prove that 
the space is an ANR, and this is often difficult.
To avoid this problem, we assume that the metric space
$Y$ is an \textit{ANRU}, which is a uniform version
of ANR. The definition of this notion is obtained from
that of the usual ANR by replacing all of ``continuous maps'' and
``neighborhoods'' by uniform ones, namely, uniformly
continuous maps and uniform neighborhoods (see \S 2.2 for
the precise definition).
This notion apparently depends on the metric of $Y$,
but in fact depends only on the uniformity it defines
(Corollary \ref{uniform invariance}). 
The assumption that $Y$ is an ANRU is not too strong,
since all compact polyhedra as well as Euclidean spaces with
the usual metric are ANRU's (Example \ref{eucl sp}). 

Now the main theorem of this paper can be stated:

\begin{Th}\label{main}
Let $X$ be a separable noncompact metrizable space and $(Y, d)$ a separable complete ANRU.
Assume that the diameters of path-components of $Y$ are bounded away from zero: that is,
the set $\{\diam C\,;\,C\text{ is a path-component of }Y\}$ has a positive lower bound.
Then, the space $C(X,Y)$ of
continuous maps from $X$ into $Y$ is an $\ell^2(2^{\aleph_0})$-manifold
with respect to the uniform topology.
\end{Th}

This theorem will be proved in \S 3.2. The condition on the path-components
of $Y$ may seem to be complicated, but this is essential as shown in \S 4.1
(Theorem \ref{components}, Corollary \ref{equiv}). 

Notice that the set $C_B(X,Y)$ of all bounded continuous maps is an open subset of $C(X,Y)$.
From examples of ANRU's mentioned above, we have the following:

\begin{Cor}\label{main cor1}
Let $X$ be a separable noncompact metrizable space and let $Y$ be
a compact polyhedron without isolated points (equipped with any admissible metric) or the
Euclidean space $\R^n\,(1\leqq n<\infty)$.
Then, $C(X,Y)$ and $C_B(X,Y)$ are $\ell^2(2^{\aleph_0})$-manifolds. 
\qed
\end{Cor}

This corollary gives many concrete examples of 
$Y$ for which $C(X,Y)$ is a Hilbert manifold for
all noncompact separable metrizable space $X$.
It can also be regarded as a ``manifold version''
of an elementary fact that the space
$C_B(X, \R)$ is a Banach space of weight $2^{\aleph_0}$,
which is homeomorphic to $\ell^2(2^{\aleph_0})$.

The next corollary is concerned with
the case where the range space is a Riemannian manifold:

\begin{Cor}\label{main cor2}
Let $X$ be a separable noncompact metrizable space, and assume that
$Y$ is a connected, complete Riemannian manifold (without boundary) equipped
with the geodesic metric. Then $C_B(X, Y)$
is an $\ell^2(2^{\aleph_0})$-manifold.
\end{Cor}

The proof of this corollary uses a variation of the main theorem
(Corollary \ref{ANRU tower}) and will be in \S 5. This corollary
would be false if one replaces $C_B(X,Y)$ by $C(X,Y)$
(Example \ref{bounded essential}). To prove our main Theorem \ref{main},
according to Toru\'{n}czyk's characterization,
we have to show that $C(X,Y)$ has a certain general position property.
This constitutes the main part of this paper (\S 3).  

The author expresses his deep gratitude to
Professors T.~Tsuboi and K.~Sakai for their invaluable help during the preparation of this paper.
He would like to thank Professor K.~Sakai also for a careful reading of the manuscript and
Professor D.~Gauld for useful comments. He also thanks K.~Mine for a useful discussion for
Appendix in this paper and S.~Saito for a careful reading of the manuscript.
Finally, he thanks the referee, who helped him a great deal to prepare
this paper in a clearer style. 

\section{Preliminaries}
\subsection{Basic terminology}
In this paper, \textit{spaces} are Hausdorff
and \textit{maps} are continuous, while
\textit{functions} are not necessarily continuous. 

Let $X$ be a space and $(Y,d)$ a metric space.
We denote by $C(X,Y)$ the space of all maps from $X$ into $Y$
with the uniform topology.

We follow the standard notation
$I=[0,1]$ and we define $\N=\{1,2,\ldots\}$.
Also, by $B(y,r)$ and $\bar{B}(y,r)$
we mean the open ball and the closed ball.
For $A\subset Y$ and $r>0$, we define $N(A, r)=\{x\in X\,;\,d(x, A)<r\}$ and
$\overline{N}(A, r)=\{x\in X\,;\,d(x, A)\leqq r\}$.

Let $(X_1, d_1), (X_2, d_2)$ be metric spaces. The \textit{product metric} on $X_1\times X_2$
is defined by $d((x_1,x_2), (x'_1, x'_2))=d_1(x_1,x'_1)+d_2(x_2,x'_2)$. 
Product spaces are assumed to have this metric and the uniformity
it induces. 

We will frequently use simplicial complexes, which are not assumed to be
finite-dimensional or locally finite. The readers are assumed to be familiar
with the basic notions concerning simplicial complexes, especially 
subdivisions. A simplicial complex is denoted by a symbol such as $K$, which is
thought of as a collection of simplexes. The union of these simplexes,
or the \textit{polyhedron} of $K$,
is denoted by $|K|$ and assumed to have the Whitehead weak topology.

We say that a function $\delta\colon (0,\infty)\to (0,\infty)$ is a \textit{modulus of continuity},
or a \textit{modulus} for short if $\delta$ is monotone and $\lim_{\varepsilon\to +0}\delta(\varepsilon)=0$.
A typical example of a modulus arises from a uniformly continuous map.
Let $f\colon X\to Y$ be a uniformly continuous map between metric spaces.
Then 
\[
\delta(\varepsilon)=\sup \{\delta\in(0,1]\,;\,d(f(x),f(x'))<\varepsilon \text{ whenever }d(x,x')<\delta\} 
\]
defines a modulus $\delta$.
We note the following:

\begin{L'a}\label{conticontimod}
For every modulus $\delta$, there exists a continuous, strictly monotone modulus $\delta'$
satisfying $\delta'(\varepsilon)\leqq \min\{1, \delta(\varepsilon)\}$ for every $\varepsilon>0$.
\end{L'a}
\begin{proof}
We can choose a strictly decreasing sequence $(\varepsilon_n)_{n\in\N}$ of positive numbers converging
to $0$ such that $\delta(\varepsilon_{n+1})<\delta(\varepsilon_n)<1$ for each $n\in\N$.
Define $\delta'\colon (0,\infty)\to (0,\infty)$ by
\[
\delta'(\varepsilon)=
\begin{cases}
(1-t)\delta(\varepsilon_{n+2})+t\delta(\varepsilon_{n+1}) \quad & \text{if }\varepsilon=(1-t)\varepsilon_{n+1}+t\varepsilon_n,\,t\in I,\\
\varphi(\varepsilon) & \text{if }\varepsilon\geqq \varepsilon_1,
\end{cases}
\]
where $\varphi$ is an arbitrary strictly monotone continuous function with
$\varphi(\varepsilon_1)=\delta(\varepsilon_2)$ and $\varphi(\varepsilon)<1$
for each $\varepsilon\geqq \varepsilon_1$.
Then $\delta'$ is continuous, strictly monotone and $\delta'\leqq \min\{1, \delta\}$.
\end{proof}

\subsection{Preliminaries on ANRU} \label{prelimANRU}
We introduce a uniform version of the notions of ANE and ANR,
namely ANEU and ANRU, but there seems to be several possible definitions for such notions
different from each other (see e.g. Isbell~\cite{Isbell} and Nguyen To Nhu~\cite{Nhu}).
We follow \cite{Nhu}, which is more convenient in our context of
metric spaces. We say that a neighborhood $V$ of a subset $A$ of a metric space
is a \textit{uniform neighborhood} if we have $N(A,\varepsilon)\subset V$ for some
$\varepsilon>0$.

\begin{Def}[ANEU]
\normalfont
A metric space $X$ is an \textit{ANEU} if, for every uniformly continuous map $f\colon A\to X$
from a closed set $A$ in a metric space $Y$, there exist a uniform neighborhood $V$
of $A$ and a uniformly continuous extension $\tilde{f}\colon V\to X$.
\end{Def}

\begin{Def}[ANRU]\label{ANEU}
\normalfont
A metric space $X$ is an \textit{ANRU} if, for every isometric closed embedding $i\colon X\to Y$
into a metric space $Y$, there exist a uniform neighborhood $V$ of $i(X)$
and  a uniformly continuous retraction from $V$ onto $i(X)$.
\end{Def}

We see that every ANRU is an ANR by a standard argument
as follows: if an ANRU $Y$ is given, we can isometrically
embed it as a closed set of a normed
space $E$ (see Bessaga-Pe\l czy\'{n}ski~\cite[Chapter II,
Corollary 1.1]{Bessaga-Pelczynski}). Since $Y$ is an ANRU,
there is an open (uniform) neighborhood $U$ of $Y$ in $E$ which
retracts onto $Y$. Since $U$ is an ANR as an open set
of a normed space, so is $Y$.

It is easy to see that every ANEU is an ANRU. A less obvious fact is that
the converse holds:

\begin{Prop}\label{ANRU=ANEU}
A metric space $(X, d)$ is an ANRU if and only if it is an ANEU.
\end{Prop}

In fact, we need only the notion of ANRU and do not need ANEU.
However, we include the proof of Proposition \ref{ANRU=ANEU}
in Appendix of this paper, since this equivalence
implies the uniform invariance of the notion of ANRU,
which is sometimes used in this paper:

\begin{Cor}\label{uniform invariance}
If $(X,d)$ and $(X',d')$ are uniformly equivalent metric spaces and $(X,d)$ is an ANRU, then so is $(X',d')$.
\end{Cor}

Indeed, it is clear from the definition that the notion of ANEU is invariant under uniform equivalences
(note that it is not so clear for ANRU). Then, this corollary follows from Propostion \ref{ANRU=ANEU}.

Now we shall give examples of ANRU's. 
A \textit{uniform complex} is the polyhedron $|K|$ of a
simplicial complex $K$, endowed with the metric
\[
d(x,y)=\sup_{v\in K^{(0)}} |x(v)-y(v)|,
\]
where $x(v)$ denotes the barycentric coordinates of $x$.
Notice that the topology induced by this metric is different from
the Whitehead weak topology unless $K$ is locally finite.
The following fact is known:

\begin{Prop}[Isbell~\cite{Isbell2}, Theorem 1.9]\label{findim uniform}
Every finite-dimensional uniform complex is an ANRU.
\qed
\end{Prop}

The next examples of ANRU's show that the main Theorem \ref{main} implies Corollary \ref{main cor1}:

\begin{Ex}\label{eucl sp}
\normalfont
The Euclidean space $\R^n$, with the usual triangulation (first subdivide the space by unit cubes and
triangulate their faces inductively, in order of increasing dimension, without introducing new vertices), is an example of a uniform complex.
As a uniform complex, $\R^n$ has a metric uniformly equivalent to
the usual one. This means, by Corollary \ref{uniform invariance} and Proposition \ref{findim uniform},
that Euclidean space $\R^n$ with the usual metric (or more generally,
a metric induced from a norm) is an ANRU.

Similarly, by Corollary \ref{uniform invariance} and Proposition \ref{findim uniform},
any compact polyhedron with a fixed admissible metric is also an ANRU,
since all admissible metrics on a compact space are uniformly equivalent.
\qed
\end{Ex}

There is an example of non-ANRU homeomorphic to $\R$:

\begin{Ex}\label{not ANRU}
\normalfont
Let $Y$ be the image of the closed embedding
$i\colon \R\to \R^2$ defined by
\[
i(x)=
\begin{cases}
(|x|, 1/x) & |x|\geqq 1, \\
(1,x)& |x|\leqq 1. \\
\end{cases}
\]   
We think of $Y$ as a metric space using the usual Euclidean metric on $\R^2$.
Then $Y$ is not an ANRU.
\qed
\end{Ex}

\subsection{Characterization of $\ell^2(\tau)$-manifolds by Toru\'{n}czyk}
Let $\tau$ be an infinite cardinal and $\Gamma$ a discrete space of cardinality $\tau$.
We denote by $\ell^2(\tau)$ the Hilbert space of weight $\tau$, in other words,
the Hilbert space whose orthonormal basis consists of $\tau$ vectors.
A metrizable space $X$ is an $\ell^2(\tau)$\textit{-manifold}
if every point of $X$ has an open neighborhood homeomorphic to
some open set in $\ell^2(\tau)$. 

The characterization of $\ell^2(\tau)$-manifolds obtained by Toru\'{n}czyk is the main ingredient of
this paper. To state this theorem, it seems to be convenient to introduce the following notion.
A map $f\colon X\to Y$ is \textit{approximable} by a map $g$ satisfying a certain property $P$
if for each map $\alpha\colon Y\to (0,\infty)$ there exists a map $g\colon X\to Y$ satisfying $P$ such that
$d(f(x),g(x))<\alpha(f(x))$ for all $x\in X$.

\begin{Th}[Toru\'{n}czyk~\cite{Torunczyk}, cf.~\cite{Torunczyk2}]
\label{characterization}
Let $\tau$ be an infinite cardinal and $\Gamma$ a discrete space of cardinality $\tau$.
The space $X$ is an $\ell^2(\tau)$-manifold if and only if it satisfies the following conditions:
\begin{itemize}
\item[(a)]$X$ is a completely metrizable ANR of weight $\tau$.
\item[(b)]For each $n\in\N$, every map from $I^n\times\Gamma$ into $X$ is approximable by
a map $g$ for which the family $\{g(I^n\times \{\gamma\})\}_{\gamma\in\Gamma}$ is discrete in $X$.
\item[(c)]For each sequence $\{K_n\}_{n\in\N}$ of finite-dimensional simplicial complexes with the cardinality of
vertices $\leqq \tau$, every map from the topological sum $\bigoplus_{n\in\N}|K_n|$ into $X$ is approximable
by a map $g$ for which the family $\{g(|K_n|)\}_{n\in\N}$ is discrete in $X$.
\end{itemize}
\end{Th}

\subsection{Basic properties of the function space}

Here we prove basic properties of the function space
$C(X,Y)$ in our main Theorem \ref{main}. 

\begin{L'a}\label{basicremark}
Let $X$ and $Y$ be separable metric spaces.
\begin{itemize}
\item[(i)] If $X$ is noncompact and $Y$ contains a homeomorphic copy of $I$,
then the weight of $C(X,Y)$ is $2^{\aleph_0}$.
\item[(ii)] If $Y$ is an ANRU then $C(X,Y)$ is an ANR.
\end{itemize}
Consequently, if $X$ and $Y$ satisfy the assumption in Theorem \ref{main},
then $C(X,Y)$ is an ANR of weight $2^{\aleph_0}$. 
\end{L'a}

\begin{proof}
(i) By the assumption, $X$ contains a countable closed discrete set $X_0$.
We think of $\partial I\subset I\subset Y$, where $\partial I$ is the
set of two endponts of $I$. Then $C(X_0, \partial I)$ is a discrete space
and its weight equals the cardinality, which is $2^{\aleph_0}$.
Consider the subset $A=\{f\in C(X,Y)\,;\,f(X_0)\subset \partial I\}$.
By the Tietze extension theorem, the restriction $A \to C(X_0,\partial I)$
is surjective. Therefore, the weight of $A$, hence of $C(X,Y)$, is at least $2^{\aleph_0}$.
The reverse inequality easily follows from the fact that both
$X$ and $Y$ are separable metrizable.

(ii) By (Bessaga-Pe\l czy\'{n}ski~\cite[Chapter II,
Corollary 1.1]{Bessaga-Pelczynski}) we can embed $Y$ isometrically
into a normed space $E$ as a closed set. Since $Y$ is an ANRU,
it has a uniform neighborhood $N$ in $E$ which has a uniformly
continuous retraction $r\colon N\to Y$. 
Then the composition of $r$ defines a (continuous) retraction $C(X,N)\to C(X,Y)$.
Notice that $C(X,N)$ is a neighborhood of $C(X,Y)$ in $C(X,E)$, which is
an AR as a normed space.
Thus $C(X,Y)$ is embedded into an AR as a retract of some neighborhood. This means
that $C(X,Y)$ is an ANR.

Finally, let $X$ and $Y$ be as in Theorem \ref{main}.
The hypothesis in (ii) is satisfied and thus $C(X,Y)$ is an ANR.
Since $Y$ is a metric space with a nontrivial path-component,
it has a homeomophic copy of $I$
(see Engelking~\cite[Excercise 6.3.12]{Engelking}).
Therefore, the hypothesis of (i) is also satisfied,
and hence $C(X,Y)$ is of weight $2^{\aleph_0}$.  
\end{proof}

\section{Proof of Theorem \ref{main}}
\label{section DAP}

In this section, we fix spaces $X$ and $Y$ as in the hypothesis of
Theorem \ref{main}: $X$ is a separable metrizable noncompact space and $(Y,d)$ is a separable
complete ANRU (\S 2.2) with the diameters of components away from zero.
For brevity, we set 
\[
\BF=C(X,Y).
\]
By Lemma \ref{basicremark}, together with the completeness of $Y$,
the function space $\BF$ satisfies the condition (a) of the
Characterization Theorem \ref{characterization} for $\tau=2^{\aleph_0}$.
The rest of this section is devoted to the proof that $\BF$ also
satisfies the conditions (b) and (c).

We may assume that $Y$ is bounded,
since if we replace $d$ by $\min\{1,d\}$,
the metric space $Y$ still has the same properties stated above and the topology of $\BF$ is
unchanged. Thus, we can define the sup-metric on $\BF$ denoted by $d$,
and the topology of $\BF$ is the one induced from $d$. 

\subsection{Preliminary lemmas for discrete approximations}
\label{preliminary lemma}

Hereafter we fix an isometric closed embedding of $(Y, d)$ into a normed space
$(E, \|\cdot\|)$, uniform neighborhood $N$ and a uniformly continuous retraction $r\colon N\to Y$
(cf. proof of Lemma \ref{basicremark}). 

Without loss of generality, we can make the following assumptions:
\begin{itemize}
\item[($\natural 1$)]Each component of $Y$ has the diameter at least $1$,
\item[($\natural 2$)]The uniform neighborhood $N$ contains the 1-neighborhood $N(Y,1)$
of $Y$ in $E$.
\end{itemize}
Indeed, these conditions are satisfied
if we multiply some positive constant to $d$.

Using the uniform continuity of $r$, we define
$\delta_1\colon(0,\infty)\to(0,\infty)$ by
\[
\delta_1(\varepsilon)=\sup\bigg\{
\delta\in(0,\varepsilon/2]\,;\,
\begin{tabular}{l}
$\|r(x)-r(x')\|<\varepsilon/2 \text{ whenever }$\\
$x,x'\in N \text{ and }\|x-x'\|<\delta$
\end{tabular}
\bigg\},
\]
Then $\delta_0$ is a modulus (see \S2.1)
satisfying $\delta_1(\varepsilon)<\varepsilon$ for all $\varepsilon>0$.
By Lemma \ref{conticontimod}, we can take a continuous, strictly monotone
modulus $\delta_0$ such that $\delta_0\leqq \min\{1, \delta_1\}$.
Consequently, we have the following:

\begin{L'a}\label{good modulus}
There exists a continuous, strictly monotone modulus
$\delta_0$ with the following properties:
\begin{itemize}
\item[(i)] $\delta_0(\varepsilon)<\min\{1,\varepsilon\}$,
\item[(ii)] $\|r(x)-r(x')\|\leqq \varepsilon/2 \text{ whenever }\|x-x'\|<\delta_0(\varepsilon),\,x,x'\in N$.
\end{itemize}
\qed
\end{L'a}

Recall that a family $\{f_\lambda\}_{\lambda\in\Lambda}$ of maps from
a space $Z$ into a metric space $(T, \rho)$
is \textit{equicontinuous} if for every $\varepsilon>0$ and $z\in Z$,
there is a neighborhood $U$ of $z$ such that
$f_\lambda(U)\subset B(f_\lambda(z), \varepsilon)$
for all $\lambda\in\Lambda$.

We introduce the following notation: 
if $F\colon Z\to \BF$ is a function,
we denote by $\check{F}\colon X\times Z\to Y$
the function defined by $\check{F}(x,t)=F(t)(x)\,(x\in X,\,t\in Z)$. 
Note that if $F$ is continuous, then so is $\check{F}$.
The converse is not true (consider the case $X=[0,\infty)$, $Y=I$,
$Z=I$, $F(t)(x)=1-(1+tx)^{-1}$).

\begin{L'a}\label{modification extension}
Let $\delta_0$ be a modulus as in Lemma \ref{good modulus}.
Let $D=\{x_n\,;\,n\in\N\}$ be a closed discrete set of $X$ and
$\{U_n\}_{n\in\N}$ a discrete open collection with $x_n\in U_n\,(n\in\N)$.
Suppose that a map $F\colon Z\to \BF$ together with
equicontinuous family $\{\gamma_n\}_{n\in\N}$ of maps $\gamma_n\colon Z\to Y$
satisfy the condition
\[ 
d(\check{F}(x_n,z), \gamma_n(z))<\delta_0
(\varepsilon(z))\,(n\in\N, z\in Z).
\]
Then for every map $\varepsilon\colon Z\to (0,1]$
there is a map $F'\colon Z\to \BF$ satisfying 
\begin{itemize}
\item[(i)] $\check{F'}(x_n, z)=\gamma_n(z)\quad(n\in\N,\,
z\in Z)$,
\item[(ii)] $d(F(z), F'(z))<\varepsilon(z)\quad(z\in Z)$,
\item[(iii)] $\check{F}|_{(X\setminus U)\times Z}=\check{F'}|_{(X\setminus U)\times Z}$,
\end{itemize}
where $U=\bigcup_{n\in\N}U_n$.  
\end{L'a}

\begin{proof}
Fix a continuous function
$\beta\colon X\to I$ with $\beta|_D=1$ and $\beta|_{X\setminus U}=0$.
Consider a function $F''\colon Z\to C(X,N)$ defined by 
\[
F''(z)(x)=
\begin{cases}
\check{F}(x,z)+\beta(x)(\gamma_n(z)-\check{F}(x_n,z)) & \text{ if } x\in U_n\subset U,\\
\check{F}(x,z)                                                    & \text{ if } x\in X\setminus U.
\end{cases}
\]
It is easy to see that $F''(z)\colon X\to E$ is certainly continuous.
To prove $F''(z)\in C(X,N)$, it is enough to check that $F''(z)(x)\in N$ for $x\in U$.
Choosing $n$ such that $x\in U_n$, we have
\[
\|\check{F}(x_n,z)-\gamma_n(z)\|<\delta_0(\varepsilon(z))<\varepsilon(z)\leqq 1,
\] 
and hence $F''(z)(x)=\check{F}(x,z)+\beta(x)(\gamma_n(z)-\check{F}(x_n,z))\in N(Y,1)\subset N$.
Thus $F''\colon Z\to C(X,N)$ is defined as a function.

We claim that $F''$ is continuous.
To see this, $z,z'\in Z$ and $x\in X$. First suppose that $x\in U_n$
for some $n\in\N$. Then we have
\[
\begin{split}
&\quad\,\,\|F''(z)(x)-F''(z')(x)\|\\
&=\|(\check{F}(x,z)+\beta(x)(\gamma_n(z)-\check{F}(x_n,z)))\\
&\quad\quad
-(\check{F}(x,z')+\beta(x)(\gamma_n(z')-\check{F}(x_n,z')))\|\\
&\leqq\|\check{F}(x,z)-\check{F}(x,z')\|+\|(\gamma_n(z)-\check{F}(x_n,z))-(\gamma_n(z')-\check{F}(x_n,z'))\|\\
&\leqq 2d(F(z), F(z'))+\|\gamma_n(z)-\gamma_n(z')\|.
\end{split}
\]
On the other hand, if $x\notin U_n$ for every $n\in\N$, we have
\[
\begin{split}
\|F''(z)(x)-F''(z')(x)\|
&=\|\check{F}(x,z)-\check{F}(x,z')\|\\
&\leqq d(F(z), F(z')).
\end{split}
\]
Combining these inequalities and
the equicontinuity of $\{\gamma_n\}_{n\in\N}$, we conclude that $F''$ is continuous.

Since $r\colon N\to Y$ is uniformly continuous, it induces
a continuous map $r_*\colon C(X,N)\to C(X,Y)=\BF$.
We define $F'\colon Z\to \BF$ by the composition
$F'=r_*\circ F''$.
Then, $F'$ is continuous and it is easy to check that $F'$
satisfies the conditions (i) and (iii).

Now it remains only to show that $F'$ satisfies (ii). 
In view of
\[
\begin{split}
\|\beta(x)(\gamma_n(z)-\check{F}(x_n, z))\|
&\leqq \|\gamma_n(z)-\check{F}(x_n, z)\|\\
&<\delta_0(\varepsilon(z)),
\end{split}
\]
we have, if $x\in U_n$,
\[
\begin{split}
&\quad\,\,\|F'(z)(x)-F(z)(x)\|\\
&=\left\|r\!\left(\check{F}(x,z)+\beta(x)(\gamma_n(z)-\check{F}(x_n, z))\right)
-r\!\left(\check{F}(x,z)\right)\right\|\\
&\leqq \varepsilon(z)/2,
\end{split}
\]
due to our choice of $\delta_0$ in Lemma \ref{good modulus}.
Therefore, in view of (iii), we have $d(F'(z),F(z))\leqq \varepsilon(z)/2<\varepsilon(z)$,
and hence, $F'$ satisfies (ii).
\end{proof}

For $\varepsilon>0$, we say that a subset $S\subset Y$ is \textit{$\varepsilon$-discrete} if for every two distinct points $p, q\in S$
we have $d(p,q)\geqq \varepsilon$, while $S$ is called \textit{$\varepsilon$-dense} if for every $y\in Y$ there exists $p\in S$
such that $d(y,p)<\varepsilon$, that is, $N(S,\varepsilon)=Y$.

Let us consider the following conditions $P_1(B,R,\varepsilon)$, $P_2(B,R,\varepsilon)$
for $B,R\subset Y$, $\varepsilon>0$:
\begin{itemize}
\item $P_1(B,R,\varepsilon):\;\; B\cap R=\emptyset\text{ and } B\cup R\text{ is }\varepsilon/4\text{\,-discrete}$.
\item $P_2(B,R,\varepsilon):\;\; \text{Both }B\text{ and }R\text{ are }\varepsilon\text{-dense}.$
\end{itemize}

We shall prove the following:

\begin{L'a}\label{blue-red extension}
Let $0<\varepsilon<1$. Then if $P_1(B,R,\varepsilon)$ holds,
there exist $B'\supset B$ and $R'\supset R$ satisfying
$P_1(B',R',\varepsilon/2)$ and $P_2(B',R',\varepsilon/2)$.
\end{L'a}

Using this lemma inductively, we have the following:

\begin{Prop}\label{blue-red points}
There exist increasing sequences $B(1)\subset B(2)\subset\cdots$
and $R(1)\subset R(2)\subset\cdots$ of 
subsets of $Y$ satisfying
$P_1(B(i), R(i), 2^{-i})$ and $P_2(B(i), R(i), 2^{-i})$ for each $i\in\N$.
\qed
\end{Prop}

It may be helpful to think of the elements of $B(i)$ and $R(i)$ as ``blue points'' and ``red points'', respectively.

\begin{proof}[Proof of Lemma \ref{blue-red extension}]
Let $B'\subset Y$ be a maximal subset satisfying
\begin{itemize}
\item $B\subset B'$,
\item $B'\cap R=\emptyset$,
\item $B'\cup R$ is $\varepsilon/8$-discrete,
\item $B'$ is $\varepsilon/4$-discrete.
\end{itemize}
(Note that $B'=B$ certainly satisfies these conditions.)
Then by the maximality,
\[
\text{for every }y\in Y,\text{ either }B(y,\varepsilon/8)\cap R\neq\emptyset\text{ or }B(y,\varepsilon/4)\cap B'\neq \emptyset.
\tag{$\sharp$}
\]
Then we have
\begin{Claim}
The set $B'$ is $\varepsilon/2$-dense. 
\end{Claim}
To show this, take $y\in Y$ with $B(y,\varepsilon/4)\cap B'=\emptyset$.
Then by ($\sharp$), we can take a point $a\in B(y,\varepsilon/8)\cap R$.
Since $\varepsilon<1$ and each path-component of $Y$ has the diameter $\geqq 1$
(see $(\natural 1)$ at the beginning of \S \ref{preliminary lemma}),
there is a point $b\in Y$ such that $d(a,b)=\varepsilon/8$.
Then $B(b, \varepsilon/8)$ does not meet $R$, because $a\in R$ and $R$ is $\varepsilon/4$-discrete.
Again by ($\sharp$), it follows that $B(b,\varepsilon/4)\cap B'\neq \emptyset$. Fix an element $c\in B(b,\varepsilon/4)\cap B'$.
Then we have
\[
d(y,c)\leqq d(y,a)+d(a,b)+d(b,c)<\varepsilon/8+\varepsilon/8+\varepsilon/4=\varepsilon/2,
\]
which proves the claim.

Let $R'\subset Y$ be a maximal subset satisfying
\begin{itemize}
\item $R\subset R'$,
\item $B'\cap R'=\emptyset$,
\item $B'\cup R'$ is $\varepsilon/8$-discrete,
\item $R'$ is $\varepsilon/4$-discrete.
\end{itemize}
(Note that $R'=R$ certainly satisfies these conditions.)
Then by an argument similar to the above, we have
\begin{Claim}
The set $R'$ is $\varepsilon/2$-dense.
\end{Claim}
Now the subsets $B', R'\subset Y$ satisfy $P_1(B',R',\varepsilon/2)$ and $P_2(B',R',\varepsilon/2)$.
\end{proof}

\subsection{Proof of the discrete approximation}
\label{discrete approximation}

Here by $\Gamma$ we denote the set of all sequences of $0$ and $1$, that is, $\Gamma=\{0,1\}^\N$.

We shall show that $\BF=C(X,Y)$ satisfies the conditions (b) and (c)
in Characterization Theorem \ref{characterization}.
These conditions will be verified at the same time. 
In fact, we prove the following:

\begin{Prop}
\label{discrete}
Let $K_\gamma\,(\gamma\in\Gamma)$ be simplicial complexes.
Then every map $\bigoplus_{\gamma\in\Gamma}f_\gamma\colon
\bigoplus_{\gamma\in\Gamma}|K_\gamma|\to \BF$ is approximable by a map $\bigoplus_{\gamma\in\Gamma}
g_\gamma$ for which $\{g_\gamma(|K_\gamma|)\}_{\gamma\in\Gamma}$
is discrete in $\BF$.
\end{Prop}

This proposition clearly implies both of the conditions (b) and (c) for $\BF$.
To prove Proposition \ref{discrete}, we show the following:

\begin{L'a}[main lemma]
\label{main lemma}
There exists a modulus $\delta$ with the following property:
for every continuous function $\alpha\colon \BF\to (0,1]$,
every family $\{K_\gamma\}_{\gamma\in\Gamma}$ of simplicial complexes
indexed by $\Gamma$ and every
map $\bigoplus_{\gamma\in\Gamma}f_\gamma
\colon\bigoplus_{\gamma\in\Gamma}|K_\gamma|\to \BF$, there exists
a map $\bigoplus_{\gamma\in\Gamma} g_\gamma\colon\bigoplus_
{\gamma\in\Gamma}|K_\gamma|\to \BF$ such that
\begin{itemize}
\item[(i)]$d(f_\gamma(p), g_\gamma(p))<\alpha f_\gamma(p)$ for every $\gamma\in\Gamma$, $p\in|K_\gamma|$ and 
\item[(ii)]$d(g_\gamma(p), g_{\gamma'}(p'))\geqq\min\{\delta\alpha f_\gamma(p), \delta\alpha f_{\gamma'}(p')\}$
for every $\gamma\neq\gamma'\in\Gamma$, $p\in|K_\gamma|$ and $p'\in|K_{\gamma'}|$.
\end{itemize}
\end{L'a}

In the above lemma, notice that there is no assumption on the simplicial complexes: they
can possess arbitrarily many vertices and can be infinite-dimensional or non-locally finite.
First we verify that Lemma \ref{main lemma} implies Proposition \ref{discrete}.

\begin{proof}[Proof of ``Lemma \ref{main lemma} $\Rightarrow$ Proposition \ref{discrete}'']
It suffices to show that condition (ii) of Lemma \ref{main lemma}
implies the discreteness of $\{g_\gamma(|K_\gamma|)\}$ in $\BF$.

If $\{g_\gamma(|K_\gamma|)\}_{\gamma\in\Gamma}$ is not discrete,
there are $F_0\in\BF$, a sequence $\{\gamma_j\}_{j\in\N}$ in $\Gamma$ with
$\gamma_j\neq \gamma_{j+1}\,(j\in\N)$, and $p_j\in |K_{\gamma_j}|\,(j\in\N)$
such that $g_{\gamma_j}(p_j)\to F_0$ as $j\to \infty$. Let $f_j=f_{\gamma_j}$ and
$g_j=g_{\gamma_j}$.
Then by (ii), for each $j\in\N$ one of the following holds:
\begin{gather*}
d(g_j(p_j), g_{j+1}(p_{j+1}))\geqq \delta\alpha f_j(p_j),\\
d(g_j(p_j), g_{j+1}(p_{j+1}))\geqq \delta\alpha f_{j+1}(p_{j+1}).
\end{gather*}
Therefore, taking a suitable subsequence, we may assume that 
$\delta\alpha f_j(p_j)\to 0$, and hence $\alpha f_j(p_j)\to 0$ (notice that $\delta$ is a modulus).
Now by (i) we have $d(f_j(p_j), g_j(p_j))\to 0$, which means $f_j(p_j)\to F_0$. Then by the continuity of $\alpha$,
we see that $\alpha f_j(p_j)\to \alpha(F_0)>0$. This contradicts $\alpha f_j(p_j)\to 0$.
\end{proof}

For a simplicial complex $K$ and its vertex $v\in K^{(0)}$, by
$\st(v,K)$ we mean the \textit{closed star} at $v$ in $K$,
that is, the union of all (closed) simplexes which has $v$ as a vertex.
Recall the following classical lemma from J. H. C. Whitehead~\cite[Theorem 35]{Whitehead2}:

\begin{L'a}[Whitehead]\label{subdivision}
For every simplicial complex $K$ and open covering $\CU$ of $|K|$,
there exists a subdivision $K'$ of $K$ such that the closed cover
$\{\st(v, K')\,;\,v\in K'^{(0)}\}$ refines $\CU$. 
\qed
\end{L'a}

\begin{proof}[Proof of Lemma \ref{main lemma}]
Fix a closed discrete subset
\[
D=\left\{ x_j^{\beta,\rho} \,;\,j\in\N,\,\beta,\rho\in\N \right\}
\]
of $X$ which is indexed by $j, \beta, \rho$. We also fix a discrete open collection
$\{U_{j}^{\beta, \rho}\}_{j,\beta,\rho}$ with $x_j^{\beta, \rho}\in U_j^{\beta, \rho}$.

We take $B(i), R(i)\subset Y\,(i\in\N)$ as in Proposition \ref{blue-red points}.
Namely, take $B(i)$ and $R(i)$ such that
\begin{itemize}
\item[(1)] $B(1)\subset B(2)\subset \cdots \subset Y,\, R(1)\subset R(2)\subset \cdots \subset Y$,
\item[(2)] $B(i)\cap R(i)=\emptyset$,
\item[(3)] $B(i)\cup R(i)$ is $2^{-(i+2)}$-discrete,
\item[(4)] Both $B(i)$ and $R(i)$ are $2^{-i}$-dense.
\end{itemize} 

Fix a modulus $\delta_0$ as in Lemma \ref{good modulus}. In particular,
$\delta_0$ satisfies the following:
\begin{itemize}
\item[(5)] $||r(x)-r(x')||\leqq \varepsilon/2$ if $x,x'\in N$ and $||x-x'||<\delta_0(\varepsilon)$,
\item[(6)] $\delta_0(\eta)<1$ for each $\eta>0$,
\item[(7)] $\delta_0$ is continuous and strictly monotone.
\end{itemize}
Then we define a modulus $\delta$ by
\begin{itemize}
\item[(8)]  $\displaystyle \delta=\frac{1}{32}\delta_0^2$,
\end{itemize}
where $\delta_0^2=\delta_0\circ \delta_0$.

Let $\alpha\colon \BF\to (0,1]$ be a continuous function,
$\{K_\gamma\}_{\gamma\in\Gamma}$ a family of simplicial complexes, and
$\bigoplus_{\gamma\in\Gamma} f_\gamma\colon \bigoplus_{\gamma\in\Gamma} |K_\gamma|\to \BF$ a map. 
We shall show that there exists a map 
$\bigoplus_{\gamma\in\Gamma} g_\gamma\colon \bigoplus_{\gamma\in\Gamma} |K_\gamma|\to \BF$
satisfying (i) and (ii) in the statement of Lemma \ref{main lemma}.

To simplify the notation, we let
\[
\alpha_\gamma=\delta_0^2\alpha f_\gamma\colon |K_\gamma|\to (0,\infty) \quad(\gamma\in\Gamma).
\]

We claim that there exists a subdivision $K'_\gamma$ of $K_\gamma$ such that
\begin{itemize}
\item[(9)] $K'_\gamma$ has the form of the second barycentric subdivision, that is,
there exists $L_\gamma$ such that $K'_\gamma=\Sd^2 L_\gamma\,\, (\gamma\in\Gamma)$,
\end{itemize}
and further that, writing $C_\gamma^v=\st(v,L_\gamma)$ for $\gamma\in\Gamma$ and $v\in L_\gamma^{(0)}$,
\begin{itemize} 
\item[(10)] $p,q\in C_\gamma^v$ implies $\displaystyle\alpha_\gamma(p)\geqq \frac{\alpha_\gamma(q)}{2}$, 
\item[(11)] $\displaystyle\diam f_\gamma(C_\gamma^v)\leqq \frac{1}{2}\inf_{p\in C_\gamma^v}\alpha_\gamma(p)$. 
\end{itemize}

To show this claim, take any $\gamma\in\Gamma$.
Let $\CU$ be the open covering of $|K_\gamma|$ defined by
\[
\CU=\left\{ \alpha_\gamma^{-1}\bigl((t,2t)\bigr)\cap f_\gamma^{-1}\!\left(B\!\left(F, \frac{\delta_0^2\alpha(F)}{8}\right)\right)
\,;\,\,t>0,\, F\in (\delta_0^2\alpha)^{-1}\bigl(\left(t,2t\right)\bigr)\right\},
\]
where $B\left(F, \frac{\delta_0^2\alpha(F)}{8}\right)$ stands for an open ball in $\BF$. It is easy to see that
the above definition certainly gives a covering of $|K_\gamma|$. By using Lemma \ref{subdivision}, we have
a subdivision $L_\gamma$ of $K_\gamma$ such that
\begin{itemize}
\item[(12)] $\{C_\gamma^v\,;\,v\in L_\gamma^{(0)}\}<\CU$,
\end{itemize}
where $C_\gamma^v=\st(v, L_\gamma)$. We shall prove that $K'_\gamma=\Sd^2 L_\gamma$
satisfies our requirements (9)--(11). Clearly, (9) is satisfied.
To see (10), let $v\in L_\gamma^{(0)}$ and $p,q\in C_\gamma^v$.
Then by (12) there exists $t>0$ such that $\alpha_\gamma(p), \alpha_\gamma(q)\in (t, 2t)$. Then
$\alpha_\gamma(q)/\alpha_\gamma(p)<2t/t=2$, which means $\alpha_\gamma(q)<2\alpha_\gamma(p)$.

To see (11), let $v\in L_\gamma^{(0)}$ and $p,q,r\in C_\gamma^v$. Then by (12) there exist $t>0$ and $F\in\BF$ such that
\begin{itemize}
\item[(13)] $\alpha_\gamma(r), \delta_0^2\alpha(F)\in (t, 2t)$ and
\item[(14)] $\displaystyle d(f_\gamma(p), F), d(f_\gamma(q), F)<\frac{\delta_0^2\alpha(F)}{8}$.
\end{itemize}
It follows from (13) and (14) that
\[
d(f_\gamma(p), f_\gamma(q))
<\frac{\delta_0^2\alpha (F)}{4}\\
<\frac{t}{2}\\
<\frac{\alpha_\gamma(r)}{2}.
\]
Since $p,q,r\in C_\gamma^v$ are arbitrary, we have (11).

\textit{Hereafter we write $K_\gamma$ for $K'_\gamma$ for simplicity.}
Let $L'_\gamma=\Sd L_\gamma$ (which means $
\Sd L'_\gamma=K_\gamma$).
Since the vertices of $L'_\gamma$ are precisely the barycenters of simplexes in $L_\gamma$,
the set $(L'_\gamma)^{(0)}$ of all vertices of $L'_\gamma$ is represented as a disjoint union
\begin{itemize}
\item[(15)] $\displaystyle (L'_\gamma)^{(0)}=\bigcup_{\theta\in\N\cup\{0\}} W_\gamma^\theta$,
\end{itemize}
where $W_\gamma^\theta$ is the set of all barycenters of $\theta$-simplexes in $L_\gamma$.
We set $D_\gamma^w=\st(w, K_\gamma)$ for $w\in (L'_\gamma)^{(0)}=\bigcup_{\theta=0}^\infty W_\gamma^\theta$.

Now define a metric $d=d_\gamma$ on the set $|K_\gamma|$ by
\begin{itemize}
\item[(16)] $\displaystyle d(p,q)=d_\gamma(p,q)=\sum_{v\in K_\gamma^{(0)}}|p(v)-q(v)|\quad(p,q\in|K_\gamma|)$,
\end{itemize} 
where $p(v)$ denotes the barycentric coordinate of $p$ at $v$.
This metric gives a topology of $|K_\gamma|$ not stronger than the original Whitehead weak topology.
We note that
\begin{itemize}
\item[(17)] if $w$ and $w'$ are distinct points in $W_\gamma^\theta$ for some
$\theta\in\N\cup\{0\}$, then $d(D_\gamma^w, D_\gamma^{w'})\geqq 2$.
\end{itemize}
As a consequence, $\{D_\gamma^w\}_{w\in W_\gamma^\theta}$ is a discrete family of closed subsets in 
$|K_\gamma|$ for each $\theta\in\N\cup\{0\}$.
By $D_\gamma(\theta)$ we denote the union of this family, that is, 
$D_\gamma(\theta)=\bigcup_{w\in W_\gamma^\theta} D_\gamma^w$. In view of (15), we observe that
\begin{itemize}
\item[(18)] $\displaystyle |K_\gamma|=
\bigcup_{w\in(L'_\gamma)^{(0)}}D_\gamma^w=\bigcup_{\theta=0}^\infty D_\gamma(\theta)$,
\end{itemize}
which is known as the dual cell decomposition of $L_\gamma$
when $L_\gamma$ is a combinatorial manifold.

As regards the relation between $D_\gamma^w \,(w\in (L'_\gamma)^{(0)})$ and $C_\gamma^v \,(v\in L_\gamma^{(0)})$,
we easily see the following:
\begin{quote}
If there exists a simplex $\sigma$ in $L_\gamma$ such that $v$ is a vertex of $\sigma$ and
$w$ is the barycenter of $\sigma$, then $\overline{N}(D_\gamma^w, 1/2)$ is contained in
$C_\gamma^v$.
\end{quote}
From this we have
\begin{itemize}
\item[(19)] $\{\overline{N}(D_\gamma^w, 1/2) \,;\, w\in (L'_\gamma)^{(0)} \}<\{C_\gamma^v \,;\,v\in L_\gamma^{(0)}\}$.
\end{itemize}

Now for each $v\in L_\gamma^{(0)}$, we let 
\begin{itemize}
\item[(20)] $\varepsilon_\gamma^v=\inf_{p\in C_\gamma^v} \alpha f_\gamma(p)$ for $v\in L_\gamma^{(0)}$,
\end{itemize}
which is positive by virtue of (10). Then define $\mu_\gamma^v\in\N$ by
\[\mu_\gamma^v=
\min\{\mu\in\N\,;\,2^{-\mu}<\delta_0^2(\varepsilon_\gamma^v)/2\}.\]  
Notice that the modulus $\delta_0$ is continuous. It is easy to see the following:
\begin{itemize}
\item[(21)] $\displaystyle \diam f_\gamma(C_\gamma^v)\leqq \frac{\delta_0^2(\varepsilon_\gamma^v)}{2}$,
\item[(22)] $\displaystyle 2^{-\mu_\gamma^v}<\frac{\delta_0^2(\varepsilon_\gamma^v)}{2}\leqq 2^{-(\mu_\gamma^v-1)}$.
\end{itemize}
For each $\theta\in\N\cup\{0\}$, we define a function $\varphi_{\gamma}^\theta\colon |K_\gamma|\to [0,1]$ 
by 
\begin{itemize}
\item[(23)] $\varphi_{\gamma}^\theta (p)=\max\{0,1-2d(p, D_\gamma(\theta))\}$.
\end{itemize}
We notice the following:
\begin{itemize}
\item[(24)] $\varphi_\gamma^\theta|_{D_\gamma(\theta)}=1$, $\varphi_\gamma^\theta|_{|K_\gamma|\setminus N(D_\gamma(\theta), 1/2)}=0$,
\item[(25)] $\varphi_{\gamma}^\theta$ is continuous with respect to metric $d$, hence to the Whitehead topology.
\end{itemize}
For each $\gamma\in\Gamma$, $j\in\N$ and $\beta, \rho\in\N$, we define a map
$h_{\gamma, j}^{\beta,\rho}\colon |K_\gamma|\to E$.
We fix $\gamma, j, \beta$ and $\rho$. For brevity, we write $\check{f}^0=\check{f}_\gamma(x_j^{\beta,\rho},\cdot)
\colon |K_\gamma|\to Y$. 

\vspace{2mm}
\noindent
\textsf{Case 1.} If $\gamma(j)=0$, we proceed as follows. In this case, we do not use $\rho$ but we use $\beta$
throughout the construction of $h_{\gamma, j}^{\beta,\rho}$.
Take any $w\in W_\gamma^\beta\subset (L'_\gamma)^{(0)}$. By (19), we have $v=v(w)\in L_\gamma^{(0)}$
such that 
\begin{itemize}
\item[(26)]$\overline{N}(D_\gamma^w, 1/2)\subset C_\gamma^{v}$.
\end{itemize}
Then by (4), we can choose $b_\gamma^w\in B(\mu_\gamma^{v})$
satisfying
\begin{itemize}
\item[(27)] $d(b_\gamma^w, \check{f}^0(C_\gamma^v))<2^{-\mu_\gamma^v}$.
\end{itemize}
Now, for every $p\in \overline{N}(D_\gamma^w, 1/2)$, using (27), (26), (21), (22) and (20), we have
\[
\begin{split}
\|b_\gamma^w-\check{f}^0(p)\|
&\leqq d(b_\gamma^w, \check{f}^0(C_\gamma^v))+\diam \check{f}^0(C_\gamma^v)\\
&<2^{-\mu_\gamma^v}+\frac{\delta_0^2(\varepsilon_\gamma^v)}{2}\\
&<\frac{\delta_0^2(\varepsilon_\gamma^v)}{2}+\frac{\delta_0^2(\varepsilon_\gamma^v)}{2}=\delta_0^2(\varepsilon_\gamma^v)\\
&\leqq \delta_0^2\alpha f_\gamma(p),
\end{split}
\]
which means that
\begin{itemize}
\item[(28)]$\|b_\gamma^w-\check{f}^0(p)\|<\alpha_\gamma(p)$ for $w\in W_\gamma^\beta$ and $p\in \overline{N}(D_\gamma^w, 1/2)$.
\end{itemize}
Notice that, by (17),
\begin{itemize}
\item[(29)] if $w\neq w'\in W_\gamma^\beta$ then $d(\overline{N}(D_\gamma^w, 1/2),\overline{N}(D_\gamma^{w'}, 1/2))\geqq 1$.
\end{itemize} 
We define $h_{\gamma, j}^{\beta, \rho}\colon |K_\gamma|\to E$ by
\begin{itemize}
\item[(30)]
$h_{\gamma, j}^{\beta, \rho}(p)=
\begin{cases}
\varphi_\gamma^\beta(p)(b_\gamma^w-\check{f}^0(p))\quad 
& \text{if }p\in \overline{N}(D_\gamma^w, 1/2) \\
& \text{ for (unique) }w\in W_\gamma^\beta, \\
& \\
0\quad & \text{if }p\in |K_\gamma|\setminus N(D_\gamma(\beta), 1/2).
\end{cases}$
\end{itemize}
By (29) and (25), we see that $h_{\gamma,j}^{\beta,\rho}$
is continuous with respect to the Whitehead topology. For 
$w\in W_\gamma^\beta$ and $p\in \overline{N}(D_\gamma^w, 1/2)$, we have
\begin{itemize}
\item[(31)]
$\|h_{\gamma, j}^{\beta, \rho}(p)\|<\alpha_\gamma(p)$
\end{itemize}
by (28). The same inequality is clearly true for $p\in |K_\gamma|\setminus \overline{N}(D_\gamma(\beta), 1/2)$.
If we define $\check{f}^0+h_{\gamma, j}^{\beta, \rho}\colon |K_\gamma|\to E$ by the
addition in $E$, we have
\begin{itemize}
\item[(32)]
$(\check{f}^0+h_{\gamma, j}^{\beta, \rho})(D_\gamma^w)=\{b_\gamma^w\}\subset B(\mu_\gamma^{v(w)})$.
\end{itemize}
by (30) and (24).

\vspace{2mm}

\noindent
\textsf{Case 2.} If $\gamma(j)=1$, we do not use $\beta$ but we use $\rho$.
For each $w\in W_\gamma^\rho$, we can take
$v=v(w)\in L_\gamma^{(0)}$ satisfying $\overline{N}(D_\gamma^w, 1/2)\subset C_\gamma^{v}$ by (19). 
By (4), we can choose
$r_\gamma^w\in R(\mu_\gamma^{v(w)})$ with $d(r_\gamma^w, \check{f}^0(C_\gamma^v))<2^{-\mu_\gamma^v}$. 
Then we can prove
\begin{itemize}
\item[(33)] $\|r_\gamma^w-\check{f}^0(p)\|<\alpha_\gamma(p)$
\end{itemize}
in the same way as (28) in Case 1.
We define a map $h_{\gamma, j}^{\beta, \rho}\colon |K_\gamma|\to E$ by
\begin{itemize}
\item[(34)]
$h_{\gamma, j}^{\beta, \rho}(p)=
\begin{cases}
\varphi_\gamma^\rho(p)(r_\gamma^w-\check{f}^0(p))\quad 
& \text{if }p\in \overline{N}(D_\gamma^w, 1/2) \\
& \text{ for (unique) }w\in W_\gamma^\rho, \\
& \\
0\quad & \text{if }p\in |K_\gamma|\setminus N(D_\gamma(\rho), 1/2).
\end{cases}$
\end{itemize}
Then by (33), for $p\in |K_\gamma|$,
\begin{itemize}
\item[(35)] $\|h_{\gamma, j}^{\beta, \rho}(p)\|<\alpha_\gamma(p)$,
\end{itemize}
which is the same inequality as (31).
If we define $\check{f}^0+h_{\gamma,j}^{\beta,\rho}\colon |K_\gamma|\to E$
using the addition in $E$,
\begin{itemize}
\item[(36)]$(\check{f}^0+h_{\gamma,j}^{\beta,\rho})(D_\gamma^w)=\{r_\gamma^w\}\subset R(\mu_\gamma^{v(w)})$ for every $w\in W_\gamma^\rho$.
\end{itemize}

\vspace{2mm}
We have defined $h_{\gamma, j}^{\beta, \rho}\colon |K_\gamma|\to E$ for every $\gamma\in\Gamma, j\in\N$
and $\beta, \rho\in \N\cup\{0\}$. Now we claim the following.
\begin{Claim}\label{equiconti}
\normalfont
$\{h_{\gamma, j}^{\beta, \rho}\}_{j,\beta,\rho}$ is equicontinuous for every $\gamma\in\Gamma$.
\end{Claim} 
First we fix $\gamma\in\Gamma$. We show the equicontinuity of
\[
\{h_{\gamma, j}^{\beta,\rho}\}_{j\in\gamma^{-1}(0),\;\beta,\rho\in\N}
\]
and that of
\[ 
\{h_{\gamma, j}^{\beta,\rho}\}_{j\in\gamma^{-1}(1),\;\beta,\rho\in\N}
\]
separately.
Since the latter family is treated similarly, we prove only the equicontunuity of the former. 
Take any point $p\in |K_\gamma|$.
It suffices to show the equicontinuity at the point $p$ with respect to the metric $d$ defined in (16).
In view of the definition (30), it suffices to consider the case $d(p,D_\gamma(\beta))\leqq 1/2$
only, since otherwise $h_{\gamma, j}^{\beta,\rho}$'s, where $j\in\gamma^{-1}(0)$,
are constantly zero in a common neighborhood of $p$.
Then there exists (unique)
$w\in W_\gamma^\beta$ such that $d(p, D_\gamma^w)\leqq 1/2$.

Let $p'\in |K_\gamma|$ be any point satisfying $d(p', p)<1/2$ and let $j\in\N$, $\beta, \rho\in\N\cup\{0\}$ and
$\check{f}^0=\check{f}_\gamma(x_j^{\beta,\rho},\cdot)\colon |K_\gamma|\to Y$.
We distinguish two cases:

\vspace{2mm}
\noindent
\textsf{Case A.} If $p'\in |K|\setminus N(D_\gamma(\beta), 1/2)$,
 then $h_{\gamma, j}^{\beta, \rho}(p')=0$ 
by the definition (30). Then by (28) and (23), we have
\[
\begin{split}
\|h_{\gamma, j}^{\beta, \rho}(p')-h_{\gamma, j}^{\beta, \rho}(p)\|
&=\|h_{\gamma, j}^{\beta, \rho}(p)\|\\
&=\varphi_\gamma^\beta(p)\|b_\gamma^w-\check{f}^0(p)\|\\
&\leqq \varphi_{\gamma}^\beta(p)\alpha_\gamma(p)\\
&=(1-2d(p, D_\gamma^w))\alpha_\gamma(p).
\end{split}
\]
Now, making use of the assumption that $d(p', D_\gamma^w)\geqq 1/2$, we have $d(p,p')\geqq 1/2-d(p, D_\gamma^w)$,
which is also expressed as $1-2d(p, D_\gamma^w)\leqq 2d(p,p')$. Combining this with the above inequality, we have
\[
\begin{split}
\|h_{\gamma, j}^{\beta, \rho}(p')-h_{\gamma, j}^{\beta, \rho}(p)\|
&\leqq 2d(p,p')\alpha_\gamma(p)=2d(p,p')\delta_0^2\alpha f_\gamma(p).
\end{split}
\]
Therefore by (6), we have
\begin{itemize}
\item[(37)] $\|h_{\gamma, j}^{\beta, \rho}(p')-h_{\gamma, j}^{\beta, \rho}(p)\|\leqq 2d(p ,p')$.
\end{itemize}

\vspace{2mm}
\noindent
\textsf{Case B.} If $p'\in \overline{N}(D_\gamma(\beta), 1/2)$, we must have 
$p'\in \overline{N}(D_\gamma^w)$ by (17), since $d(p', p)<1/2$ and $d(p, D_\gamma^w)<1/2$.  
Consider $v(w)\in L_\gamma^{(0)}$
and $b_\gamma^w\in B(\mu_{\gamma}^{v(w)})$
that are chosen when we defined $h_{\gamma, j}^{\beta, \rho}$. 
Then we have
\[
\begin{split}
\|h_{\gamma, j}^{\beta, \rho}(p')-h_{\gamma, j}^{\beta, \rho}(p)\|
&=\|\varphi_\gamma^\beta(p)(b_\gamma^w-\check{f}^0(p))-\varphi_{\gamma}^\beta(p')(b_\gamma^w-\check{f}^0(p'))\|\\
&\leqq|\varphi_{\gamma}^\beta(p)-\varphi_\gamma^\beta(p')|\,\|b_\gamma^w-\check{f}^0(p)\|\\
&\qquad+|\varphi_\gamma^\beta(p')|\,\|\check{f}^0
(p)-\check{f}^0(p')\|\\
&= 2|d(p, D_\gamma^w)-d(p', D_\gamma^w)|\,\|b_\gamma^w-\check{f}^0(p)\|\\
&\qquad +|\varphi_\gamma^\beta(p')|\,\|\check{f}^0(p)-\check{f}^0(p')\|\\ 
&\leqq 2d(p,p')\|b_\gamma^w-\check{f}^0(p)\|+d(f_\gamma(p), f_\gamma(p'))\\
&\leqq 2d(p,p')\alpha_\gamma(p)+d(f_\gamma(p), f_\gamma(p')),
\end{split}
\]
where the last inequality follows from (28). By (6), we have
\begin{itemize}
\item[(38)]
$\|h_{\gamma, j}^{\beta, \rho}(p')-h_{\gamma, j}^{\beta, \rho}(p)\|
\leqq 2d(p,p')+d(f_\gamma(p), f_\gamma(p'))$.
\end{itemize}
The inequalities (37) and (38) shows that the family
$\{h_{\gamma,j}^{\beta,\rho}\}_{j\in\gamma^{-1}(0),\,\beta,\rho\in\N}$ is equicontinuous at 
$p\in |K_\gamma|$.
The proof of Claim \ref{equiconti} is complete.

\vspace{3mm}
We define $g_{\gamma, j}^{\beta, \rho}\colon |K_\gamma|\to Y
\,(\gamma\in\Gamma,\,j\in\N,\,\beta,\rho\in\N)$ by
\begin{itemize}
\item[(39)]$g_{\gamma,j}^{\beta,\rho}(p)=r(\check{f}^0(p)+h_{\gamma,j}^{\beta,\rho}(p))$,
\end{itemize}
where $\check{f}^0=\check{f}_\gamma(x_j^{\beta,\rho},\cdot)\colon |K_\gamma|\to Y$ and 
the addition is performed in $E$. By (35) and (6) we have $\|h_{\gamma, j}^{\beta, \rho}(p)\|<\alpha_\gamma(p)
<1$, hence 
\[
\check{f}^0(p)+h_{\gamma,j}^{\beta,\rho}(p)\in N(Y,1)\quad(p\in |K_\gamma|),
\]
which shows that $g_{\gamma,j}^{\beta,\rho}$ is well-defined. Since $r$ is a uniformly continuous
retraction onto $Y$, we conclude from Claim \ref{equiconti} that
\begin{itemize}
\item[(40)]
$\{g_{\gamma, j}^{\beta, \rho}\}_{j,\beta,\rho}$ is equicontinuous for every $\gamma\in\Gamma$.
\end{itemize}
Notice that $\check{f}_\gamma(x_j^{\beta,\rho},p)=\check{f}^0(p)=r(\check{f}^0(p))\,\,(p\in |K_\gamma|)$. 
Then by (31), (35) and (5),
\begin{itemize}
\item[(41)]
$d(g_{\gamma, j}^{\beta,\rho}(p), \check{f}_\gamma(x_j^{\beta,\rho},p))
<\delta_0\alpha f_\gamma(p)$.
\end{itemize}
Furthermore, by (32), (36) and the definition (39 of $g_{\gamma, j}^{\beta, \rho}$, we have
\begin{itemize}
\item[(42)] $g_{\gamma, j}^{\beta, \rho}(D_\gamma^w)=\{b_\gamma^w\}\subset B(\mu_\gamma^{v(w)})$
if $\gamma(j)=0$ and $w\in W_\gamma^\beta$,
\item[(43)] $g_{\gamma, j}^{\beta, \rho}(D_\gamma^w)=\{r_\gamma^w\}\subset R(\mu_\gamma^{v(w)})$
if $\gamma(j)=1$ and $w\in W_\gamma^\rho$.
\end{itemize}

By Lemma \ref{modification extension}, (40) and (41), for each $\gamma\in\Gamma$ there exists 
a map $g_\gamma\colon |K_\gamma|\to \BF$ satisfying
\begin{itemize}
\item[(44)] $\check{g}_\gamma(x_j^{\beta,\rho},\cdot)=g_{\gamma, j}^{\beta, \rho}\quad(j\in\N,\,\beta,\rho\in \N)$,
\item[(45)] $d(g_\gamma(p), f_\gamma(p))<\alpha f_\gamma(p)\quad(p\in |K_\gamma|)$,
\end{itemize}
and $\check{g}_\gamma|_{(X\setminus U)\times |K_\gamma|}=\check{f}_\gamma|_{(X\setminus U)\times |K_\gamma|}$, where $U=\bigcup_{j,\beta,\rho}U_j^{\beta,\rho}$.
The above inequality (45) means that the condition (i) is satisfied.
In the rest of the proof, we show that the condition (ii) is also satisfied. 

Let $\gamma\neq\gamma'\in\Gamma$, $p\in |K_\gamma|$ and $p'\in |K_{\gamma'}|$. We may assume that
$\gamma(j_0)=0$ and $\gamma'(j_0)=1$ for some $j_0\in\N$. By (18, there exist $\beta_0, \rho_0\in
\N$, $w_0\in W_\gamma^{\beta_0}$ and $w'_0\in W_\gamma^{\rho_0}$ such that $p\in D_\gamma^{w_0}$
and $p'\in D_{\gamma'}^{w'_0}$. 

Let $x_0=x_{j_0}^{\beta_0, \rho_0}$,
$g_0=g_{\gamma, j_0}^{\beta_0,\rho_0}$ and $g'_0=g_{\gamma', j_0}^{\beta_0, \rho_0}$. 
Moreover we let $D_0=D_\gamma^{w_0}$ and $D'_0=D_{\gamma'}^{w'_0}$.
Then by (44)
we have $\check{g}_\gamma(x_0,\cdot)=g_0$ and $\check{g}_{\gamma'}(x_0,\cdot)=g'_0$.
Furthermore by (42) and (43)
\begin{itemize}
\item[(46)]
$g_0(p)=b_\gamma^{w_0}\in B(\mu_\gamma^{v(w_0)})$,
\end{itemize}
since $\gamma(j_0)=0$. Similarly, by (43)2,
\begin{itemize}
\item[(47)]
$g'_0(p')=r_{\gamma'}^{w'_0}\in R(\mu_\gamma^{v(w'_0)})$,
\end{itemize}
since $\gamma'(j_0)=1$. We let $v_0=v(w_0)\in L_\gamma^{(0)}$ and $v'_0=v(w'_0)\in L_{\gamma'}^{(0)}$.
Further, we let $\mu_0=\mu_\gamma^{v_0}$ and $\mu'_0=\mu_{\gamma'}^{v'_0}$.
By (1), (2) and (3) we have
\begin{itemize}
\item[(48)] $\displaystyle d(g_0(p), g'_0(p'))\geqq d(B(\mu_0), R(\mu'_0))\geqq 
\frac{1}{4}\min\{2^{-\mu_0}, 2^{-\mu'_0}\}$.
\end{itemize}
On the other hand, by (10), we have
\begin{itemize}
\item[(49)] $\alpha_\gamma(q)\geqq \alpha_\gamma(p)/2$ for every $q\in C_\gamma^{v_0}$ and
\item[(50)] $\alpha_{\gamma'}(q')\geqq \alpha_{\gamma'}(p')/2$ for every  $q'\in C_{\gamma'}^{v'_0}$.
\end{itemize}
Then, by using (46)--(50), (22), (20) as well as (7)
(to exchange $\delta_0$ and infimum), we have the following:
\[
\begin{split}
d(g_\gamma(p), g_{\gamma'}(p'))
&\geqq d(\check{g}_\gamma(x_0, p), \check{g}_{\gamma'}(x_0, p'))=d(g_0(p), g'_0(p'))\\
&\geqq \frac{1}{4}\min\{2^{-\mu_0}, 2^{-\mu'_0}\}\\
&\geqq \frac{1}{4}\cdot\frac{1}{4}\{\delta_0^2(\varepsilon_\gamma^{v_0}), \delta_0^2(\varepsilon_{\gamma'}^{v'_0})\}\\
&=\frac{1}{16}\min\left\{\inf_{q\in D_0}\alpha_{\gamma}(q), \inf_{q'\in D'_0} \alpha_{\gamma'}(q')\right\}\\
&\geqq \frac{1}{32}\min\{\alpha_\gamma(p), \alpha_\gamma(p')\}\\
&=\frac{1}{32}\min\{\delta_0^2\alpha f_\gamma(p), \delta_0^2\alpha f_{\gamma'}(p')\}\\
\end{split}
\]
By (8), this inequality means 
\[d(g_\gamma(p), g_{\gamma'}(p'))\geqq 
\min\{\delta\alpha f_\gamma(p), \delta\alpha f_{\gamma'}(p')\}.
\]
This show that the condition (ii) is
satisfied. The proof is complete.  
\end{proof}

This lemma implies Proposition \ref{discrete}, as is seen before.
Then it follows that the function space $\BF$ satisfies the conditions (b) and (c). 
We have verified that $\BF$ satisfies the conditions (a), (b), (c) in Theorem \ref{characterization}
and the proof of Theorem \ref{main} is completed.

\section{Remarks on Theorem \ref{main}}
\subsection{The role of the condition on the components}
In Theorem \ref{main}, the assumption
that $Y$ has the diameters of components away from zero
is essential in the following sense: 

\begin{Th} \label{components}
Let $(Y,d)$ be a separable complete ANRU. If $C(\N, Y)$ is an $\ell^2(2^{\aleph_0})$-manifold,
then the diameters of components of $Y$ are bounded away from zero.
\end{Th}

\begin{proof}
First notice that $Y$ has no isolated points. Indeed, if $Y$ has
a isolated point $p$, the constant map $\N\to Y$ with the value $p$
is an isolated point in $C(\N,Y)$, which means that $C(\N,Y)$
is not an $\ell^2(2^{\aleph_0})$-manifold, contradicting to the assumption. 

Thus, to prove our conclusion, it is enough to consider the case where $Y$ has
infinitely many path-components $C_i\,(i\in\N)$.

First we show that $\{d(C_i, C_j)\,;\,i\neq j\in\N\}$ has a positive lower bound.
Suppose this is not true. We can embed $Y$ isometrically into
a normed space $E$ as a linearly independent closed set (Bessaga-Pe\l czy\'{n}ski~\cite[Chapter II,
Corollary 1.1]{Bessaga-Pelczynski}). For $i<j$,
take a line segment $L_{i,j}$ in $E$ connecting a point in $C_i$ and a
point in $C_j$ such that $\diam L_{i,j}<2d(C_i, C_j)$.
Then the diameters of $L_{i,j}$'s do not have positive
lower bound.
Since $Y$ is linearly independent in $E$, the segments
$L_{i,j}$ are disjoint except for their endpoints: that is,
if $p\in L_{i,j}\cap L_{k,l}$ for $(i,j)\neq (k,l)$,
then $p$ is a common endpoint of $L_{i,j}$ and $L_{k,l}$.
Moreover, each $L_{i,j}$ intersects $Y$ only in their endpoints.

Consider the subset $Z=Y\cup \bigcup_{i<j} L_{i,j}$ of $E$.
Then $Y$ is closed in $Z$. Since $Y$ is an ANRU, there exist
a uniform neighborhood $N$ and a retraction $r\colon N\to Y$. 
Then for some $i,j$, the segment $L_{i,j}$ is contained in $N$. Then,
the image $r(L_{i,j})$ is contained in $Y$
and intersects both $C_i$ and $C_j$. This means that
$r(L_{i,j})$ is not path-connected, a contradiction.

Thus, the set $\{d(C_i, C_j)\,;\,i\neq j\in\N\}$
has a lower bound $\delta>0$. To obtain the conclusion,
suppose that the diameters of components are not bounded
away from zero. Then there exists a subsequence $\{C'_i\}_{i\in\N}$
of $\{C_i\}_{i\in\N}$ such that $\diam C'_i\to 0$ as $i\to\infty$.
Hereafter we simply write $C_i$ to mean $C'_i$.
For each $i\in\N$,
fix a point $p_i\in C_i$ and a countable dense subset
$D_i$ of $C_i$.
Consider an element $f_0\in C(\N,Y)$ defined by $f_0(i)=p_i$
and a subset 
\[
U=\{f\in C(\N,Y)\,;\,f(i)\in C_i\text{ for all }i\in\N\}
\]
of $C(\N,Y)$. Then $U$ is a neighborhood of $f_0$ in
$C(\N,Y)$. Indeed, $U$ contains the $\delta$-neighborhood
of $f_0$. Moreover, $U$ is separable. Indeed,
\[
\bigcup_{j\in\N}\{f\in C(\N,Y)\,;\,f(i)\in D_i\text{ if }i\leqq j\text{ and }f(i)=p_i\text{ if }i>j\}
\]
is a countable dense subset of $U$, since $\diam C_i\to 0$ as $i\in\infty$.
This means that $f_0\in C(\N,Y)$ has a separable neighborhood,
which contradicts the assumption
that $C(\N,Y)$ is an $\ell^2(2^{\aleph_0})$-manifold.
\end{proof}

Thus, we have a refinement of Theorem \ref{main} as follows:

\begin{Cor}\label{equiv}
Let $(Y,d)$ be a separable complete ANRU. The following are equivalent:
\begin{itemize}
\item[(1)] $C(X,Y)$ is an $\ell^2(2^{\aleph_0})$-manifold for every
separable noncompact metrizable space $X$.
\item[(2)] $C(\N, Y)$ is an $\ell^2(2^{\aleph_0})$-manifold.
\item[(3)] The diameters of path-components of $Y$ are bounded away from zero.
\end{itemize}
\qed
\end{Cor}

\subsection{Remarks on the space of bounded maps and counterexamples}\label{conclusion}
Consider the subspace $C_B(X,Y)$ of all bounded continuous
maps from $X$ to $Y$. Then $C_B(X,Y)$ is open in $C(X,Y)$,
and hence $C_B(X,Y)$ is an $\ell^2(2^{\aleph_0})$-manifold if $C(X,Y)$ is.
The converse is false, and this will be shown in Example \ref{kuchibashi}. 
The following will be used in the proof of
Corollary \ref{main cor2}:

\begin{Cor}\label{ANRU tower}
Let $X$ be a separable noncompact metrizable space and $(Y,d)$ a metric space.
Suppose that $Y$ can be represented as $Y=\bigcup_{i\in\N} Y_i$, where 
$\{Y_i\}_{i\in\N}$ satisfies the following properties:
\begin{itemize}
\item $\{Y_i\}_{i\in\N}$ is an increasing sequence of complete ANRU's,
\item For each $i\in\N$, the diameters of the path-components of $Y_i$ are bounded away from zero,
\item For each $i\in\N$, we have either $d(Y_i, Y\setminus Y_{i+1})>0$ or $Y_i=Y$.
\item Every bounded subset of $Y$ is contained in some $Y_i$.
\end{itemize}
Then $C_B(X,Y)$ is an $\ell^2(2^{\aleph_0})$-manifold.
\end{Cor}

\begin{proof}
We can think of $C_B(X,Y_i)$ as a subset of $C_B(X,Y)$.
By the first two assumptions for $\{Y_i\}$ and Theorem \ref{main}, $C_B(X,Y_i)$ is an $\ell^2(2^{\aleph_0})$-manifold
for each $i\in\N$. The third assumption shows that $C_B(X,Y_{i+1})$ is a neighborhood of
$C_B(X,Y_i)$ in $C_B(X,Y)$. The last assumption means $C_B(X,Y)=\bigcup_{i\in\N} C_B(X,Y_i)$.
Thus, we have the conclusion. 
\end{proof}

\begin{Ex}\label{kuchibashi}
\normalfont
Let $Y\subset \R^2$ be the non-ANRU homeomphic to $\R$
as in Example \ref{not ANRU}. If we let
$Y_n=\{(x,y)\in Y\,;\,x\leqq n\}$, then $\{Y_n\}_{n\in\N}$
satisfies the conditions in Corollary \ref{ANRU tower}.
Therefore for every separable noncompact
metrizable space $X$, the space $C_B(X,Y)$ of bounded continuous
functions is an $\ell^2(2^{\aleph_0})$-manifold.

However, the space $C(X,Y)$ is not necessarily a $\ell^2(2^{\aleph_0})$-manifold.
Indeed, $C(\N,Y)$ is not locally connected at $f_0$ defined by $f_0(m)=(m,1/m)$.
To see this, take any $\varepsilon > 0$ and choose $n \in \N$ so that $2/n < \varepsilon$.
Then define $f \in C(\N,Y)$ as follows:
$f(n) = (n, -1/n)$ and $f|_{\N \setminus \{n\}} = f_0|_{\N\setminus \{n\}}$.
Then, $d(f_0,f) = 2/n < \varepsilon$ but
$f_0$ and $f$ are not connected by a path with diameter $<1$.
\qed
\end{Ex}

The compact-open topology is another important topology on $C(X,Y)$.
We write $C_k(X,Y)$ to denote $C(X,Y)$ with this topology, while $C(X,Y)$
means the space with the uniform topology as before.

\begin{Ex}\label{nonmanifold}
\normalfont
The space $C_k(X,Y)$ can fail to be locally path-connected (and hence to be a Hilbert manifold), 
even if $C(X,Y)$ is a Hilbert manifold.

For example, let
\[
X=Y=\bigcup_{n\in\Z}\{(x,y)\in\R^2\,;\,(x-2n)^2+y^2=1\}\subset \R^2.
\]
Then, $Y$ can be regarded as a metric space as a subset of $\R^2$
and is uniformly equivalent to a 1-dimensional
uniform complex (see Proposition \ref{findim uniform}).
This means that $Y$ is an ANRU by Corollary \ref{uniform invariance}.
Thus, $C(X,Y)$ is a $\ell^2(2^{\aleph_0})$-manifold by Theorem \ref{main}.
 
We shall prove that $C_k(X,Y)$ is not locally path-connected at the identity map $X\to Y$.
For every neighborhood $V$ of the identity in $C_k(X,Y)$, there exists $n\in\N$ such that
the map $c_n\colon X\to Y$ belongs to $V$, where $c_n$ is defined by
\[
c_n(x,y)=
\begin{cases}
(x,y)\quad &\text{if $|x|\leqq n$},\\
(n,0)\quad &\text{if $x\geqq n$},\\
(-n,0)\quad &\text{if $x\leqq -n$}.
\end{cases}
\]
Then $c_n$ is not homotopic to the identity, that is, there is no path in $C_k(X,Y)$
from $c_n$ to the identity.
\qed
\end{Ex}

\begin{Ex}
\normalfont
The space $C(X,Y)$ can fail to be locally path-connected
even if $C_k(X,Y)$ is homeomorphic to the Hilbert space. 
Let $X, Y$ be as in Example \ref{kuchibashi}. Then, $C_k(X,Y)\approx
C_k(\N,\R)\approx \R^\N\approx \ell^2$, where the last homeomorphism
is proved in Anderson-Bing~\cite{Anderson}, whereas $C(X,Y)$ is
not locally path-connected as seen before.
\qed
\end{Ex}

A popular notion similar to ANRU is that of \textit{uniform ANR}
introduced by Michael~\cite{Michael}. Every ANRU is clearly
a uniform ANR. We do not know if our assumption that ``$Y$ is an ANRU''
can be loosened to ``$Y$ is a uniform ANR'', so we ask the next
question:

\begin{Q}
\normalfont
Let $X$ be a separable, noncompact metric space and let $(Y,d)$
be a separable, path-connected complete metric space that is a uniform ANR.
Then, is $C(X,Y)$ an ANR? Is it an $\ell^2(2^{\aleph_0})$-manifold?
\end{Q}

\section{Bounded maps into complete Riemannian manifolds}
In this section we apply Corollary \ref{ANRU tower} to the case where
$Y$ is a connected complete Riemannian manifold (of class $C^\infty$)
to prove Corollary \ref{main cor2}.  
Let $K$ be a nonempty compact subset of $M$. We say that $K$ is \textit{strictly convex} if
every two distinct points $p,q\in K$ can be joined by a unique geodesic contained in $K$ and this geodesic
is contained in the interior $\Int K$ except for its endpoints.   
Clearly if two strictly convex sets meet, then the intersection
is strictly convex.
The next lemma is known in the folklore and is easy to prove:

\begin{L'a}\label{strictly convex with nonempty interior}
Let $M$ be an $n$-dimensional Riemannian manifold without boundary.
If $V\subset M$ is strictly convex and has nonempty interior in $M$, then
$V$ is homeomorphic to $I^n$.
\qed
\end{L'a}

We use the following gluing theorem for ANRU's:

\begin{Th}[Nguyen To Nhu~\cite{Nhu2}]\label{pasting}
Let $(X,d)$ be a metric space and $A_1, A_2$ closed subspaces with $A_1\cap A_2\neq\emptyset$ and $A_1\cup A_2=X$.
Then $(X,d)$ is an ANRU if the following (i) and \textrm{(ii)} hold:
\begin{itemize}
\item[i)] $A_1,A_2$ and $A_1\cap A_2$ are ANRU's, 
\item[ii)] the metric $d$ on $X$ is uniformly equivalent to $d'$, where $d'$ is defined by
\[
d'(x,y)=
\begin{cases}
d(x,y) \quad & x,y\in A_i\,(i=1,2),\\
\inf\{d(x,z)+d(y,z)\,;\,z\in A_1\cap A_2\} \quad & \text{otherwise}.
\end{cases}
\]
\end{itemize}
\end{Th}

Notice that if $X$ is compact, the condition (ii) for $d'$ is always satisfied.
In this case, we can even drop the condition $A_1\cap A_2\neq \emptyset$.

\begin{proof}[Proof of Corollary \ref{main cor2}]
Let $M$ be a connected complete Riemannian manifold. It suffices to show that $M$
satisfies the assumption for $Y$ in Corollary \ref{ANRU tower}.
We can cover $M$ with locally finite collection of
countably many strictly convex sets $V_i\,(i\in\N)$ such that
the interiors $U_i=\Int V_i$ are nonempty.
Indeed, we can take such $V_i$'s to be closed balls $V_i=\bar{B}(p_i,r_i)$,
with the interior $U_i=\Int V_i$ being the open ball $B(p_i,r_i)$
(see Lang~\cite[Chapter VIII, Theorem 5.8 and Lemma 5.9]{Lang}).
Further we assume that
$V_{i_1}\cap\cdots \cap V_{i_n}\neq \emptyset$ implies $U_{i_1}\cap\cdots \cap U_{i_n}\neq \emptyset$
for any finitely many $i_1,\ldots,i_n$. This is realized by shrinking each $V_i$ slightly.

Let $M_i=\bigcup_{k=1}^i V_k\,(i\in\N)$. 
We show that every $M_i$ is an ANRU by induction on $i$.
For $i=1$, $M_1=V_1$ is strictly convex set with nonempty interior
and hence is homeomorphic to $I^n$ by
Lemma \ref{strictly convex with nonempty interior}.
Thus $V_1$ is an ANRU by Proposition \ref{findim uniform}.
For $i=2$, $V_1$ and $V_2$ are ANRU's by the case $i=1$.
If $V_1$ does not meet $V_2$, then $M_2=V_1\cup V_2$ is an ANRU as noticed above.
If $V_1$ meets $V_2$,  then $U_1$ also meets $U_2$,
and hence $V_1\cap V_2$ is again a strictly convex set with nonempty interior, which is
an ANRU. By Theorem \ref{pasting}, $M_1=V_1\cup V_2$ is an ANRU.
Suppose that $i\geqq 3$ and the assertion is true for $i-1$. Then $M_i=M_{i-1}\cup V_i$ and  
\[
M_{i-1}\cap V_i=\bigcup_{k=1}^{i-1} (V_k\cap V_i). 
\]
In the above, the right hand side is an ANRU by the inductive hypothesis.
Then it follows that $M_j$ is an ANRU by Theorem \ref{pasting}. 

We can replace $\{M_i\}$ with a suitable subsequence such that
$M_i\subset \Int M_{i+1}$ for each $i$.
Then, $\{M_i\}$ satisfies the conditions for $\{Y_i\}_{i\in\N}$ in Corollary \ref{ANRU tower}.
Indeed, the last condition holds since any bounded
closed set in a complete Riemannian manifold is compact 
(see Lang~\cite[Chapter VIII, Corollary 6.7]{Lang}).
\end{proof}

\begin{Ex}\label{bounded essential}
\normalfont
In Corollary \ref{main cor2}, we cannot replace $C_B(X,Y)$ by $C(X,Y)$.
We can modify the embedding $i\colon\R\to \R^2$ given in Example \ref{not ANRU} slightly
to get a $C^\infty$-embedding $i'\colon \R\to \R^2$
which agrees with $i$ outside small neighborhoods of $1$ and $-1$. Then, the image $Y=i'(\R)$ is
naturally a complete connected Riemannian manifold.
The same argument as in Example \ref{kuchibashi} shows that the space
$C(\N,Y)$ is not an $\ell^2(2^{\aleph_0})$-manifold.
\qed
\end{Ex}

\section{Generalization of the domain space}
Neither the separability nor the metrizability of the domain space
$X$ is essential,\footnote{This remark is due to D. Gauld.}
but still some restrictions corresponding to noncompactness of $X$,
which is necessary to take a closed discrete set of ``sufficiently many''
points in $X$, should be made. 
We can also loosen the condition on the density of $Y$ without essential change.

Recall that the \textit{density} of a space is the minimal cardinality of dense sets and that
the \textit{extent} of a space is the supremum of the cardinalities of closed discrete sets. For metrizable
spaces they both coincide with the weight of the space. We say that a space \textit{attains} its extent
if there exists a closed discrete subset whose cardinality equals to the extent. Note that a separable noncompact
metric space attains its extent $\aleph_0$. 

\begin{L'a}\label{weightgen}
Let $X$ be a normal space whose density $\lambda$ is equal to the extent and $Y$ be a metrizable space
of weight $\nu$ which has a nontrivial path-component. If $X$ attains its extent and $\nu\leqq\lambda$,
then the weight of $C_B(X,Y)$ is equal to $2^\lambda$.
\end{L'a} 

\begin{proof}
Analogous to Lemma \ref{basicremark} (i).
Notice that if $\nu\leqq\lambda$ then $2^\lambda\leqq\nu^\lambda\leqq 
(2^\nu)^\lambda=2^{\nu\lambda}=2^\lambda$.
\end{proof}

By this lemma, Theorem \ref{main} is generalized as follows:

\begin{Th}
Let $X$ be a normal space whose density $\lambda$ is equal to the extent
and $(Y, d)$ an complete ANRU of weight not greater than $\lambda$.
If $X$ attains its extent and the diameters of path-components of
$Y$ are bounded away from zero, then
$C_B(X,Y)$ is an $\ell^2(2^{\lambda})$-manifold.
\qed
\end{Th}

If the domain space $X$ is metrizable,
the property that $X$ attains its extent is described by a covering property similar to
compactness. According to Barbati and Costantini~\cite{Barbati-Costantini},
a metrizable space $Z$ is \textit{generalized compact} (abbriviated GK) if
every open cover of $Z$ has a subcover of cardinality (strictly) less than its weight.  
Costantini~\cite[Proposition 2.3]{Costantini} has proved that a
metrizeble space attains its extent if and only if it is not GK.
Hence, we have the following:

\begin{Th}
Let $X$ be a metrizable, not GK space of weight $\lambda$
and $(Y, d)$ an complete ANRU of weight not greater than $\lambda$.
If the diameters of path-components of $Y$ are bounded away from zero,
then $C_B(X,Y)$ is an $\ell^2(2^{\lambda})$-manifold.
\qed
\end{Th}

\section*{Appendix. The equivalence of ANRU and ANEU}
\renewcommand{\thesection}{A}
\setcounter{Th}{0}
In this appendix, we prove Proposition \ref{ANRU=ANEU} which states
the equivalence of ANRU and ANEU.\footnote{The
following argument is essentially due to Nguyen To Nhu~\cite{Nhu}.
The author is also indebted to K.~Mine for useful suggestions for this appendix.}
The ``if'' part of the proposition is easy. To prove the ``only if'' part, it suffices to show the following:

\begin{L'a}\label{custom-made extensor}
Let $X$, $Y$ be metric spaces and $A\subset Y$ a nonempty closed set.
For every uniformly continuous map $f\colon A\to X$,
there exist an isometric closed embedding $X\hooklongrightarrow Z$ into a metric space $Z$ and
a uniformly continuous extension $\tilde{f}\colon N(A,\varepsilon)\to Z$ for some $\varepsilon>0$. 
\end{L'a}

The next lemma is a special case when the attaching map is an isometry. 

\begin{L'a}[Nguyen To Nhu~\cite{Nhu2}]\label{adjunction space}
Let $(X, d_X)$, $(Y,d_Y)$ be metric spaces and $A\subset X$ a nonempty closed set.
Let $f\colon A\to Y$ be an isometric embedding.
Consider the disjoint union $X\sqcup Y$ as a set and the equivalence relation $\sim$
generated by $a\sim f(a)$, $a\in A$. Then we can define a metric $d$ on the set
$Z=X\sqcup Y/\sim$ as follows:
\[
d(x,y)=
\begin{cases}
d_X(x,y) \quad & x,y\in X, \\
d_Y(x,y)         & x,y\in Y, \\
\inf_{a\in A} (d_X(x,a)+d_Y(f(a),y)) \quad & x\in X, y\in Y.
\end{cases}
\] 
Furthermore, $Y$ is naturally isometrically embedded as a closed set in $Z$
and $f$ naturally extends to an isometric embedding $\tilde{f}\colon X\to Z$.
\qed
\end{L'a}

We prepare the results required to show Lemma \ref{custom-made extensor}.
First, we need a basic result on ANEU's, which is proved by Isbell~\cite[Theorem 1.4]{Isbell}.\footnote{It
should be noticed that Isbell~\cite{Isbell} uses the term ANEU in the setting of
uniform spaces, and Isbell's ANEU is more general when restricted to metric spaces.}

\begin{L'a}\label{real line is ANEU}
The real line $\R$ is an ANEU with respect to the usual metric.
\qed
\end{L'a}

Let $D$ be a set and $Y$ a metric space. Define $F(D,Y)$ to be the space of all functions
from $D$ into $Y$ endowed with a metric $d(f,g)=
\min\{1, \sup_{x\in D} d(f(x), g(x))\}$. 
We think of $D$ as a metric space by the standard discrete metric.
Since a fucntion $f\colon Z\to F(D,Y)$ from a metric space $Z$
is uniformly continuous if and only if 
the corresponding function $\check{f}\colon D\times Z\to Y$
is uniformly continuous, we can easily show the following:

\begin{L'a}\label{functions into ANEU}
If $Y$ is an ANEU, then $F(D,Y)$ is an ANEU.
In particular, $F(D, \R)$ is an ANEU. 
\qed
\end{L'a}

This lemma enables us to prove:

\begin{Prop}\label{pseudometric extension}
Let $(X, d)$ be an metric space and $\rho$ a uniformly continuous
pseudometric on a closed subset $A\subset X$.
Then, for some $\varepsilon>0$, the pseudometric
$\rho$ can be extended to a uniformly continuous
pseudometric $\tilde{\rho}$ on $N=N(A,\varepsilon)$ 
such that $A$ is closed in $(N, \tilde{\rho})$.  
\end{Prop}

\begin{proof}
Define $g\colon A\to F(A,\R)$ by $g(x)(y)=\rho(x,y)$.
Since $\rho$ is uniformly continuous, $g$ is
uniformly continuous. By Lemma \ref{functions into ANEU},
there exist $\varepsilon>0$ and a uniformly continuous extension
$\tilde{g} \colon N=N(A,\varepsilon)\to F(A,\R)$.
Notice that, for $x,y\in A$, 
\[
\rho(x,y)=\sup_{a\in A}|\tilde{g}(x)(a)-\tilde{g}(y)(a)|.
\]
Define $\tilde{\rho}_0\colon N\times N\to \R$ by the same formula above
(this is well-defined for sufficiently small $\varepsilon$).
By definition $\tilde{\rho}_0$ is an extension of $\rho$ and a uniformly continuous pseudometric on $N$.
With respect to this pseudometric, the subset $A$ may not be closed.
We can define the desired pseudometric $\tilde{\rho}$ by 
$\tilde{\rho}(x,y)=\tilde{\rho}_0(x,y)+|d(x, A)-d(y, A)|$.
\end{proof}

Now we can prove Lemma \ref{custom-made extensor} to finish the proof
of Proposition \ref{ANRU=ANEU}.

\begin{proof}[Proof of Lemma \ref{custom-made extensor}]
Define $\rho\colon A\times A\to [0,\infty)$ by $\rho(x,y)=d(f(x),$ $f(y))$. 
Then $\rho$ is a uniformly continuous pseudometric on $A$.
By Proposition \ref{pseudometric extension}, for some $\varepsilon>0$, we can extend $\rho$ to
a uniformly continuous pseudometric $\tilde{\rho}$ on $N=N(A,\varepsilon)$ such that $A$ is closed
in $(N, \tilde{\rho})$. Consider the quotient metric spaces $E=N/\tilde{\rho}$ and $B=A/\rho$.
We can consider $B$ as a closed subset of $E$.
The function $g\colon B\to X$, induced by $f$, is an isometric embedding.
Then, by the Lemma \ref{adjunction space}, we can embed $X$ isometrically
into $Z=E\sqcup X/\!\sim$ as a closed subset and $g$ extends to an isometric embedding
$\tilde{g}\colon E\to Z$. The desired extension $\tilde{f}\colon N\to Z$
can be defined as the composition $N\twoheadrightarrow E\stackrel{\tilde{g}}{\rightarrow} Z$.
\end{proof}

\end{document}